\documentclass[dblspace,leqno]{aopcm}
\newcommand{\ink}{\rule{.2\baselineskip}{.25\baselineskip}}
\newcommand{\noi}{\noindent }

\newcommand{\la}{\lambda  }

\newcommand{\ti}{\tilde }
\newcommand{\cc}{\cdot }
\newcommand{\non}{\nonumber }
\newcommand{\dist}{{\rm dist}}

\newcommand{\ii}{\iota }
\newcommand{\In}{\mbox{In} }
\newcommand{\pa}{\partial }
\newcommand{\ba}{\begin{array}}
\newcommand{\ea}{\end{array}}
\newcommand{\bea}{\begin{eqnarray}}
\newcommand{\eea}{\end{eqnarray}}
\newcommand{\beas}{\begin{eqnarray*}}
\newcommand{\eeas}{\end{eqnarray*}}
\newcommand{\be}{\begin{equation}}
\newcommand{\ee}{\end{equation}}
\newcommand{\bc}{\begin{center}}
\newcommand{\ec}{\end{center}}
\newcommand{\ben}{\begin{enumerate}}
\newcommand{\een}{\end{enumerate}}

\newcommand{\ei}{\end{itemize}}

\newcommand{\ds}{\displaystyle}

\newcommand{\brmk}{\begin{remark}\begin{em}}
\newcommand{\ermk}{\end{em}\end{remark}}
\newcommand{\bexa}{\begin{example}\per\begin{em}}
\newcommand{\eexa}{\end{em}\end{example}}

\newcommand{\skp}{\vspace{\baselineskip}}

\newcommand{\E}{I\!\!E}

\newcommand{\acal}{{\cal {A}}}

\newcommand{\I}{{\cal {I}}}

\newcommand{\F}{{\cal {F}}}
\newcommand{\C}{{\cal {C}}}
\newcommand{\G}{{\cal {G}}}
\newcommand{\D}{\Delta}

\newcommand{\Om}{\Omega}
\newcommand{\om}{\omega}

\newcommand{\del}{\delta}

\newcommand{\R}{I\! \! R}
\newcommand{\INT}{I\! \! N}
\newcommand{\ov}{\overline}
\newcommand{\Df}{\doteq}

\newcommand{\half}{\frac{1}{2}}
\newcommand{\beq}{\begin{eqnarray*}}
\newcommand{\eeq}{\end{eqnarray*}}
\newcommand{\beqn}{\begin{eqnarray}}
\newcommand{\eeqn}{\end{eqnarray}}
\newcommand{\inn}[2]{\langle {#1}, {#2}\rangle}

\newcounter{theorem}
\newcounter{remark}
\newcounter{lemma}
\newcounter{cor}
\newcounter{defn}
\newcounter{cond}
\newcounter{assu}
\newcounter{prop}
\newcounter{pr}
\newcounter{definition}

\newtheorem{theorem}{Theorem}[section]
\newtheorem{remark}{Remark}[section]
\newtheorem{lemma}{Lemma}[section]

\newtheorem{defn}{Definition}[section]

\newtheorem{assu}{Condition}[section]

\begin{document}

\SPECFNSYMBOL{}{1}{2}{}{}{}{}{}{}%
\AOPMAKETITLE

\AOPyr{2000}
\AOPvol{00}
\AOPno{00}
 \AOPpp{000--000}
\AOPReceived{Received }
 \AOPAMS{Primary 60J60; secondary 60J65, 60K25, 34D20}
 \AOPKeywords{Stability,
positive recurrence, invariant measures, Skorokhod problem,
constrained processes, constrained ordinary differential
equation, queueing systems, law of large numbers.}
\AOPtitle{On Positive Recurrence of Constrained Diffusion Processes}%
\AOPauthor{Rami Atar, Amarjit Budhiraja\thanks{Research supported in part by the National
Science Foundation (NSF-DMI-9812857) and the University of Notre
Dame Faculty Research Program} and Paul Dupuis\thanks{Research supported in part
by the National Science Foundation (NSF-DMS-9704426 and DMS-0072004) and the Army Research Office (ARO-DAAD19-99-1-0223)}}%
\AOPaffil{Technion-Israel Institute of Technology, University of North Carolina at Chapel Hill and Brown University}%
\AOPlrh{RAMI ATAR, AMARJIT BUDHIRAJA AND PAUL DUPUIS}%
\AOPrrh{STABILITY OF CONSTRAINED DIFFUSIONS}%
\AOPAbstract{Let $G \subset \R^k$ be a convex polyhedral cone
with vertex at the origin given as the intersection of half
spaces $\{G_i, i= 1, \cdots , N\}$, where $n_i$ and $d_i$ denote
the inward normal and direction of constraint associated with
$G_i$, respectively.  Stability properties of a class of
diffusion processes, constrained to take values in $G$, are
studied under the assumption that the Skorokhod problem defined
by the data $\{(n_i, d_i), i = 1, \cdots , N \}$ is well posed
and the Skorokhod map is Lipschitz continuous. Explicit
conditions on the drift coefficient, $b(\cdot)$, of the diffusion
process are given under which  the constrained process is
positive recurrent and has a unique invariant measure.  Define
$$ \C \Df \left \{ - \sum_{i=1}^N \alpha_i d_i; \alpha_i \ge 0, i
\in \{1, \cdots , N\} \right \}. $$ Then the key condition for
stability is that there exists $\delta \in (0, \infty)$ and a
bounded subset $A$ of $G$ such that for all $x \in G\backslash
A$, $b(x) \in \C$ and $\dist(b(x), \partial \C) \ge \delta $,
where $\partial \C$ denotes the boundary of $\C$.}%

\maketitle



\section{Introduction}%
The stability properties of constrained stochastic processes are
of central importance in the study of queuing systems that arise
in computer networks, communications and manufacturing problems.
In recent years there has been a significant progress in the study
of stability of such systems \cite{harwil, harwil2, dupwil,
cheman, dai, dai2, che, mal, meyn, BD}. All the papers in the
list above which treat the heavy traffic diffusion model consider
the case where both the drift and the diffusion coefficients are
constant. However, in many applications a homogeneous model of this
kind is not well suited, and it is enough if we mention systems where there is a control that depends on the system's state.
Although a similar motivation leads one
to also study variable constraint directions on the boundary, we confine
ourselves here to fixed directions. In fact, some basic
stability properties of the corresponding Skorokhod map, that we take
advantage of here, are not yet well understood in the setting of variable directions of constraint.

In this paper we consider the stability properties of constrained
diffusion processes when both the drift and the diffusion
coefficients may be state dependent. Let $G \subset \R^k$ be a
convex polyhedral cone with vertex at the origin given as the
intersection of half spaces $\{G_i, i= 1, \cdots , N\}$, where
$n_i$ and $d_i$ denote the inward normal and direction of
constraint associated with $G_i$ respectively.  The stochastic
processes considered in this paper will be constrained to take
values in $G$.  One of our central assumptions is that the
Skorokhod problem defined by the data $\{(d_i, n_i); i = 1, \cdots
, N \}$ is well posed on all of $D_G([0, \infty): \R^k)$ (the
space of functions $\phi$ which are right continuous, have left
limits and $\phi(0) \in G$) and the Skorokhod map $\Gamma : D_G([0,
\infty): \R^k) \to D_G([0, \infty): \R^k)$ is Lipschitz
continuous. We refer the reader to \cite{harrei, dupish1, dupram1}
for sufficient conditions under which the Skorokhod map is
Lipschitz continuous.
The paper \cite{rei} studies networks of single class queues for which the Skorokhod problem associated to a diffusion approximation is regular.
Some examples of feedforward networks which
lead to a regular Skorokhod problem have been studied in \cite{ngu,
pet}.  An example of a multiclass networks with feedback which leads to
a regular Skorokhod problem has recently been studied in
\cite{dupram3}.

In this paper we consider the constrained diffusion process
$\{X^x(t)\}_{t \ge 0}$
 given as the unique solution of the
equation: \be \label{maineq} X^x(t) = \Gamma\left (x +
\int_0^{\cdot} \sigma(X^x(s)) dW(s) + \int_0^{\cdot}b(X^x(s)) ds
\right) (t); \;\; t \ge 0 \ee We assume global Lipschitz
conditions on $\sigma$ and $b$ (cf.\ (\ref{lipbsig})), the boundedness
(\ref{growsig}) and uniform non degeneracy (Condition \ref{irred})
 of the diffusion coefficient $\sigma$.  The main
result we present is that under the above conditions the Markov
process $\{X^x(t)\}_{t \ge 0}$ is positive recurrent and has a
unique invariant measure if there exists a $\delta  \in (0,
\infty)$ and a bounded subset $A$ of $G$ such that for all $x \in
G \backslash A$, $b(x) \in \C(\delta)$, where
$$\C(\del) \Df \{v \in \C:{\rm dist}(v,\partial C)\ge\del\},$$
and
$$
\C \Df \left \{ - \sum_{i=1}^N \alpha_i d_i; \alpha_i \ge 0,
i \in \{1, \cdots , N\} \right \}
$$
and $\partial \C$ denotes the boundary of $\C$.
The
Lipschitz conditions on $\sigma$ and $b$  are assumed to assure that there is a unique
solution
to the constrained diffusion process (\ref{maineq}).  As we point
out in Remark \ref{finrem}, our main result continues to hold if
these Lipschitz and growth assumptions are replaced by the
assumption
that (\ref{maineq}) has a unique weak solution for every $x \in
G$
and the solution is a Feller Markov process.
It should also be observed that the non-degeneracy assumption on $\sigma$ and the Feller property are used only in Section 4 in proving the ergodicity of the
constrained
diffusion process,
and are not needed to prove stability.

As in \cite{dupwil,mal,dai} the key idea in the proof of this
result is to study stability properties of a related deterministic
dynamical system.  In Theorem \ref{staba} we show that for each $\delta > 0$, the family of
deterministic constrained trajectories defined as
\be
z(t) \Df \Gamma \left (x + \int_0^{\cdot} v(s) ds \right )(t),
\label{dettyp} \ee
for which $v(t) \in \C(\delta)$, $t \in [0, \infty)$,  enjoys strong uniform stability properties
if for some $\delta \in (0, \infty)$, $v(t) \in \C(\delta)$
for all $t \in  [0, \infty)$. These stability
properties enable us to use $T(x)$, the hitting time to the
origin (cf. (\ref{hittime}), as a  Lyapunov function for the
stability analysis of the stochastic problem. We study some basic
properties of $T(\cc)$ in Lemma \ref{hit}. The key consequence of
the stability of (\ref{dettyp}) is Lemma
\ref{hit} (iii).  This result along with the Lipschitz property of
the Skorokhod map leads to Lemma \ref{deltastep} which is the
crucial step in relating the stability of the deterministic
dynamical system (\ref{dettyp}) with that of the stochastic system
(\ref{maineq}).
Stability and instability results for reflecting Brownian motion have been obtained in a number of different settings.
The papers \cite{wil,hobrog} consider the two dimensional case with constant and non-constant directions of constraint on the two faces of the domain,
respectively.
In \cite{harwil2} conditions are presented which guarantee the stability of a multi-dimensional reflecting Brownian motion,
and in addition characterizes conditions under which the invariant distribution  has a product form distribution.
The paper \cite{che} also obtains sufficient conditions for stability for this class of processes.

Our main result is Theorem \ref{finitethm} where
it is shown that there is a compact set $B \in G$ for
which the hitting time: $$ \tau_{B}(x) \Df \inf \{t:
X^x(t) \in B\}$$ has finite expectation which as a
function of $x$ is bounded on compact subsets of  $G$.
Proofs of such results generally use a Lyapunov function that is in the domain of the generator of the process (e.g., twice continuously differentiable in our case).
An interesting feature of the approach we use is that far less regularity is required of the Lyapunov function.
The paper concludes with the proof of the positive recurrence and
the uniqueness of the invariant measure for $\{X^x(t)\}$, which is
standard due to the uniform non-degeneracy of the diffusion
coefficient.

\section{Definitions and Formulation}
Let $G \subset \R^k$ be the convex polyhedral cone in $\R^k$ with
the vertex at origin given as the intersection of half spaces
$G_i$, $i = 1, \cdots , N$. Let $n_i$ be the unit vector
associated with $G_i$ via the relation
$$ G_i = \{ x \in \R^k :
\inn{x}{n_i} \ge 0 \}.$$
Denote the boundary of a set $S \subset
\R^k$ by $\partial S$. We will denote the set $\{ x \in \pa G:
\inn{x}{n_i} = 0\}$ by $F_i$. For $x \in \pa G$, define the set,
$n(x)$, of unit inward normals to $G$ at $x$ by $$ n(x) \Df \{r: |r| =
1, \;\;\; \inn{r}{x-y} \le 0, \;\; \forall y \in G\}.$$ With each
face $F_i$ we associate a  unit
vector $d_i$ such that $\inn{d_i}{n_i} > 0$.
This vector defines the {\it direction of constraint} associated with the face $F_i$.
For $x \in \pa
G$ define
$$ d(x) \Df \left\{d \in \R^k: d = \sum_{i \in \In(x)}\alpha_i d_i;
\alpha_i \ge 0; \;\; |d| = 1 \right\}, $$ where $$\In(x) \Df \{i \in
\{1,2, \cdots N\}: \inn{x}{n_i} = 0 \}.$$ Let $D([0,\infty):\R^k)$
denote the set of functions mapping $[0,\infty)$ to $\R^k$ that
are right continuous and have limits from the left.  We endow
$D([0,\infty):\R^k)$ with the usual Skorokhod topology. Let
\[
 D_G([0,\infty):\R^k) \doteq \{ \psi \in
D([0,\infty):\R^k) : \psi(0) \in G\}.
\]
For $\eta \in D([0,\infty):\R^k)$
let $| \eta | (T)$
denote the total variation of $\eta$ on
$\left[ 0,T \right]$ with respect to the Euclidean norm on
$\R^k$.

\begin{defn}
\label{def-sp}
Let $\psi \in D_G([0,\infty):\R^k)$
be given.
Then $(\phi, \eta) \in D([0,\infty):\R^k)\times D([0,\infty):\R^k)$
 solves the Skorokhod problem (SP) for $\psi$
with respect to $G$ and $d$  if and only if $\phi(0) = \psi (0)$,
 and for all
$t \in [0, \infty)$
\begin{enumerate}
\item
$\phi(t) = \psi (t) + \eta (t)$;
\item
$\phi (t) \in G$;
\item
$| \eta | (t) < \infty$;
\item
$\ds | \eta | (t) = \int_{[0,t]}
 I_{\left\{ \phi (s) \in \partial G \right\} }
d | \eta | (s)$;
\item
There exists Borel measurable $\gamma : [ 0, \infty)
\rightarrow \R^{k}$
such that $\gamma (t) \in d( \phi (t))$,
$d|\eta |$-almost everywhere
and
\[ \eta (t) = \int_{[0,t]}
\gamma (s)
d | \eta | (s). \]
\end{enumerate}
\end{defn}
On the domain $D \subset  D_G([0,\infty):\R^k)$ on which there is a
unique solution to the Skorokhod problem we define the Skorokhod
map (SM) $\Gamma$ as $\Gamma(\psi) \Df \phi$, if $(\phi, \psi - \phi)$
 is the unique
solution of the Skorokhod problem posed by $\psi$.
We will make the following assumption on the regularity
of the Skorokhod map defined by the data $\{(d_i, n_i); i = 1, 2,
\cdots N\}$.
\begin{assu}
\label{regular}
The Skorokhod map is well defined on all of
$D_G([0,\infty):\R^k)$, i.e.,
$D = D_G([0,\infty):\R^k)$ and the SM
is Lipschitz continuous in the following sense.
There exists a $K < \infty$ such
that
for all $\phi_1, \phi_2 \in D_G([0,\infty):\R^k)$:
\be
\sup_{0 \le t < \infty}|\Gamma(\phi_1)(t) - \Gamma(\phi_2)(t)|
< K \sup_{0 \le t < \infty}|\phi_1(t) - \phi_2(t)|.
\label{lip}
\ee
\end{assu}
We will assume without loss of generality that $K \ge 1$.
We refer the reader to \cite{dupish1} (or alternatively see  \cite{dupram1})
for sufficient conditions under which this regularity property holds.

We now introduce the
constrained diffusion process that will be
studied in this paper. Let $(\Om, \F , P)$ be a complete
probability space on which is given a filtration $\{\F_t\}_{t \ge
0}$ satisfying the usual hypotheses.  Let $(W(t), \F_t)$ be a
$n$-dimensional standard Wiener process on the above probability
space. We will study the constrained diffusion process
given as a solution to
  equation \ref{maineq}, namely,
\[
\label{constrainz}
X^x(t) = \Gamma\left (x + \int_0^{\cdot} \sigma(X^x(s)) dW(s) +
\int_0^{\cdot}b(X^x(s))
ds \right) (t),\]
where $\sigma : G \to \R^{k\times k}$ and
$b: G \to \R^k$ are maps satisfying the following condition:
\begin{assu}
\label{growlip}
There exists $\gamma \in (0, \infty)$ for which
\be
\label{lipbsig}
|\sigma(x) - \sigma(y)| + |b(x) - b(y)| \le \gamma|x-y|;
\;\;\; \forall x,y \in G
\ee
and
\be
\label{growsig}
|\sigma(x)| \le \gamma; \;\;\;
\forall x \in G .
\ee
\end{assu}
Using the regularity assumption on the Skorokhod map it can be
shown (cf.\ \cite{andore,dupish1}) that there is a well defined
process satisfying (\ref{constrainz}).
In fact,
the classical method of Picard iteration gives the following:
\begin{theorem}
\label{uniqsoln}
For all $x \in G$ there exists a unique pair of continuous $\{\F_t\}$ adapted
processes
$(X^x(t), k(t))_{t \ge 0}$
and a
progressively measurable process $(\gamma(t))_{t \ge 0}$ such that
the following hold:
\begin{enumerate}
\item
$X^x(t) \in G$, for all $t \ge 0$, a.s.
\item  For all $t \ge 0$,
$$X^x(t) = x + \int_0^t \sigma(X^x(s)) dW(s) + \int_0^t b(X^x(s)) ds +
k(t),$$
a.s.
\item
For all $T \in [0, \infty)$
$$
|k|(T) < \infty , \;\; a.s.
$$
\item
$$|k|(t) = \int_0^t I_{\{X^x(s) \in \pa G\}} d|k|(s),$$
and
$k(t) = \int_0^t \gamma(s) d|k|(s)$ with $\gamma(s) \in
d(X^x(s))$
a.e. $[d|k|]$.
\end{enumerate}
\end{theorem}
\begin{remark}\label{rem:sm}
The process $X^x(\cdot)$ is the unique continuous $\{\F_t\}$ adapted process
which satisfies the equation
$$
X^x(t) = \Gamma\left (x + \int_0^{\cdot} \sigma(X^x(s)) dW(s) +
\int_0^{\cdot}b(X^x(s))
ds \right) (t),$$
for all $t$ a.s.  Also, $X^x(\cc)$ is a Feller Markov process.
\end{remark}
We now proceed to formulate our central result.
Define
\[
\C \Df \left \{ - \sum_{i=1}^N \alpha_i d_i: \alpha_i \ge 0;
\; i \in \{1, \cdots , N\} \right \}
\]
The cone ${\cal C}$ was used to characterize stability of certain
semimartingale reflecting Brownian motions in \cite{BD}.

For $\delta \in (0, \infty)$, define
$$
\C(\delta) \Df
\{v \in \C:\dist(v,\partial C)\ge\del\}.
$$
Our next assumption on the diffusion model stipulates the
permissible velocity directions.
\begin{assu}
\label{permit}
There exist a $\delta  \in (0, \infty)$  and a bounded set
$A \subset G$ such that
for all $x \in G\backslash A$, $b(x) \in \C(\delta)$.
\end{assu}
Finally we will make the following uniform  nondegeneracy assumption on
the diffusion coefficient.
\begin{assu}
\label{irred}
There exists $c \in (0, \infty)$ such that
for all $x \in G$ and $\alpha \in \R^k$
$$
\alpha'(\sigma(x)\sigma'(x))\alpha \ge c \alpha'\alpha .
$$
\end{assu}

Here is the main theorem of this paper.
\begin{theorem}
Assume that Conditions \ref{regular}, \ref{growlip}, \ref{permit},
\ref{irred}
hold.
Then the strong Markov process $\{X^x(\cc); x \in G\}$ is
positive recurrent and has a unique invariant probability
measure.
\label{main}
\end{theorem}
In the rest of the paper we will assume that Conditions
\ref{regular}, \ref{growlip}, \ref{permit} and
\ref{irred}
hold.

\section{Stability of Constrained ODEs}
Let $v: [0, \infty) \to \R^k$ be a measurable map such that
\be
\int_0^t |v(s)| ds < \infty; \;\;\; \mbox{for all} \;\;
t \in [0, \infty).
\label{boc}
\ee
Let $x \in G$.  In this section we will study the stability
properties of the trajectory $z: [0, \infty) \to \R^k$ defined as
\be
\label{traj}
z(t) \Df \Gamma \left(x + \int_0^{\cdot} v(s) ds\right)(t); \;\; t \in [0,
\infty).
\ee

It is useful to rewrite the above trajectory  as a solution
of
an ordinary differential equation.  In order to do so we
introduce the following notion of discrete projections
(cf. \cite{
dupish1,dupram1}).
Define $\pi: \R^k \to G$ as follows:
$$
\pi(y) \Df \Gamma(\psi_y)(1), \;\;\; y \in \R^k,$$
where $\psi_y \in D([0, \infty); \R^k)$ is given as
\[
\psi_y(t) = \left\{ \begin{array}{l} 0, \;\;\; t \in [0,1) \\
y, \;\;\; t \in [1, \infty). \end{array} \right.
\]
In other words,
$\pi$ is a projection that is consistent with the given Skorokhod problem,
in that the constrained version of any piecewise constant trajectory $\psi$ can be found by recursively applying $\pi$.
We also define
 the {\it projection of the velocity} $v\in \R^k$ {\it at} $x\in G$ by
$$
\pi(x,v) \Df \lim_{\Delta \to 0} \frac{\pi(x+\Delta v) -
  x}{\Delta}.
$$
For a proof of the fact that the above limit exists we refer the
reader to \cite{BD} where various properties of the
projection map are also studied.  In particular we will use the
following facts (for proofs see \cite{BD}).
\begin{enumerate}
\item
For $x \in G$,
$\alpha, \beta , \gamma \ge 0$ and $v \in \R^k$:
\be
\pi(\beta x , \alpha v + \gamma x) = \alpha \pi(x, v) + \gamma x.
\label{scal}
\ee
\item
For $v \in \R^k$, we have that
\be
\pi(v) = 0 \;\; \mbox{if and only if} \;\;
v  \in \C
\label{stayz}
\ee
\end{enumerate}

The following theorem represents the trajectory in (\ref{traj})
as a solution of an ordinary differential equation
(cf. \cite{dupish1}).
\begin{theorem}
\label{ode2}
Let $v: [0, \infty) \to \R^k$ satisfy (\ref{boc}).  Then for all
$x \in G$,
$z(\cdot)$ defined via (\ref{traj}) is the
unique absolutely continuous function
such
that
\[
\dot{z}(t) = \pi(z(t), v(t)), \;\;\; a.e. \;\; t, \;\;\;\;
z(0) = x.
\]
\end{theorem}

In Theorem \ref{staba} below we present a basic stability property of the
above dynamical system.
\begin{theorem}
\label{staba}
Let $v$ be as in Theorem \ref{ode2}.
Assume that there exists a $\delta  \in (0, \infty)$
such that
$$
v(t) \in \C(\delta) \;\; \mbox{for all}\;\;
 t \in [0,
\infty).
$$
Let $x \in G$ and $z(\cdot)$ be defined via (\ref{traj}).
Then:
$$
|z(t)| \le \frac{K^2|x|^2}{K|x| + \delta t}, \;\;\; \forall
t \in [0, \infty),
$$
where $K$ is the finite constant in (\ref{lip}).
\end{theorem}

\noi
{\bf Proof:}
In order to specify the initial point of the trajectory we will
write
the trajectory defined by (\ref{traj}) as $z(x, \cc)$.
Define the trajectory $\ti z(\cc)$ as
$$\ti z(t) \Df \Gamma\left( \int_0^{\cc}v(s) ds\right)(t).$$
Theorem \ref{ode2} implies that $\ti z(\cc)$ is the unique
solution
of
\be
 \dot{\ti z}(t) = \pi(\ti z(t), v(t)), \;\; a.e.\; t; \;\;\;\;
\ti z(0) = 0.
\label{subs}
\ee
However
since $v(t) \in \C(\delta) \subset \C$,
we have  from (\ref{stayz})
that $\pi(0, v(t))  =0$ for all $t \ge 0$ and
so the zero trajectory solves (\ref{subs}).
By Theorem \ref{ode2} this implies that
$\ti z(t) \equiv 0$.
Thus
\beqn
\sup_{0 \le t < \infty}
|z(x,t)| &=&
\sup_{0 \le t < \infty}
|z(x,t) - \ti z(t)| \non
\\
&=&
\sup_{0 \le t < \infty}
\left|\Gamma\left(x + \int_0^{\cc}v(s) ds\right)(t)
- \Gamma \left(\int_0^{\cc}v(s) ds\right )(t)\right| \non \\
&\le &
K |x|.
\label{zero}
\eeqn

The above inequalities in particular show that the theorem is
true
when $x = 0$. Henceforth we assume that $x \neq 0$.
Define
$$
\gamma \Df \frac{\delta }{K|x|}$$
and
\be
\label{psidef}
\psi(t) \Df (1+ \gamma t)z(x,t).\ee
Note that from (\ref{scal}) it follows
that
\beq
\dot{\psi}(t) &=& \gamma z(x,t) + (1+ \gamma t)\pi(z(x,t),
v(t))\\
&=&
\pi\left(z(x,t), (1+\gamma t) v(t) + \gamma z(x,t)\right) \\
&=&
\pi\left(\psi(t), (1+\gamma t)\left(v(t) + \frac{\gamma}{1 +
  \gamma t} z(x,t)\right)\right).
\eeq
By Theorem \ref{ode2} we now have that
$\psi(t) = \Gamma(x + f(\cc))(t)$ for all $t \in [0, \infty)$
where
$$
f(t) \Df \int_0^t\left( (1+\gamma s)\left(v(s) + \frac{\gamma}{1 +
  \gamma s} z(x,s)\right)\right) ds.
$$
Note that from (\ref{zero}) it follows that
$$
\left|\frac{\gamma}{1 +
  \gamma t} z(x,t)\right| \le K|x| \frac{\delta }{K|x|} = \delta .$$
Thus
if $v \in \C(\delta )$
then
$v  + \frac{\gamma}{(1 +
  \gamma t)} z(x,t) \in \C$.  From this observation it follows
that
for all $t \in [0, \infty)$
$$
u(t) \Df (1+\gamma t)\left(v(t) + \frac{\gamma}{1 +
  \gamma t} z(x,t)\right) \in \C.$$

Define the trajectory
$$\ti \psi(t) \Df \Gamma (f(\cc))(t); \;\; t \in [0, \infty).$$
Then $\ti \psi(\cc)$ solves the equation
\be
\label{dummy}
\dot{\ti \psi}(t) = \pi(\ti \psi(t), u(t)); \;\; \ti \psi(0) = 0.
\ee
Since for all $t \in (0, \infty)$, $u(t) \in \C$
we have that $\pi(0, u(t)) = 0$.  Thus the function
$x(t) = 0$ for all $t \in (0, \infty)$ is a solution of
(\ref{dummy}).  Now by the uniqueness of the solution of
(\ref{dummy}) (Theorem \ref{ode2}) we have that
$\ti \psi(t) = 0$ for all $t \in (0, \infty)$.
Thus
\beq
|\psi(t)| &=& |\psi(t) - \ti \psi(t)|\\
& \le & |\Gamma(x + f(\cc))(t) - \Gamma(f(\cc))(t)|\\
& \le & K|x|,
\eeq
for all $t \in (0, \infty)$.
Finally from (\ref{psidef})
$$
|z(x,t)| \le \frac{K|x|}{1+ \gamma t} = \frac{K^2 |x|^2}{
K|x| + \delta  t}.
$$
\ink

For $x \in G$ and $\delta \in (0, \infty)$
let $\acal(x, \delta )$ be the collection of all absolutely
continuous functions $z: [0, \infty) \to \R^k$
defined via (\ref{traj}) for some $v: [0, \infty) \to
{\cal C}(\del)$
which satisfies (\ref{boc}).
Henceforth we will fix such a $\delta $ and abbreviate
$\acal(x, \delta )$ by $\acal(x)$.

For a fixed $x \in G$, we now define the
``hitting time to the origin'' function
as follows:
\be
\label{hittime}
T(x) \Df \sup_{z \in \acal(x)} \inf \{t \in [0, \infty): z(t) =
0\}.\ee
We next study some of the properties of $T(x)$.
\begin{lemma}
\label{hit}
There exist constants $c,C\in(0,\infty)$ depending only on
$K$ and $\del$ such that the following holds.
\begin{itemize}
\item[(i)]
For all $x, y \in G$
$$
|T(x) - T(y)|  \le C|x-y|.
$$
\item[(ii)]
\[
 T(x) \ge c|x|. \]
Thus, in particular, for all $M \in (0, \infty)$
the set $\{x \in G: T(x) \le M \}$ is compact.
\item[(iii)] Fix $x \in G$ and
let $z \in \acal(x)$.  Then for all $t > 0$
\[
T(z(t)) \le (T(x) - t)^+.
\]
\end{itemize}
\end{lemma}
{\bf Proof:}
We first show that for all $x \in G$
\be
\label{bdd}
T(x) \le \frac{4K^2}{\delta} |x|.\ee
Fix $x \in G$ and let $z \in \acal(x)$ be arbitrary.  From Theorem \ref{staba}
we have that for all $t \in (0, \infty)$
$$
|z(t)| \le \frac{K^2|x|^2}{K|x| + \delta t}.
$$
Hence for all $t\ge T_1\doteq 2K^2\delta^{-1}|x|$ one has
$|z(t)| \le |x|/2$.
In general, if $$T_n \Df T_1\sum_{k=0}^{n-1}2^{-k}$$ then for $t\ge T_n$
one has that $$|z(t)| \le \frac{|x|}{2^n}.$$ Thus $z(t) = 0$
for all $t \ge 4K^2\delta^{-1}|x|$.
Since $z \in \acal(x)$ is arbitrary, (\ref{bdd}) follows.

Now let $x, y \in G$ be arbitrary.  Let
$\{z_n\} \subset \acal(x)$ be a sequence such that if
$$
\tau_n \Df \inf \{t: z_n(t) = 0 \},$$
then
\be
\label{limit}
\tau_n \to T(x) \;\;\; \mbox{as}\;\; n \to \infty .
\ee
Note that $z_n$ is given as
$$z_n(t) = \Gamma \left(x + \int_0^{\cc}v_n(s) ds\right)(t);
\;\;\;
t \in [0, \infty)
$$
for some $v_n$ satisfying (\ref{traj}).
Define for $t \in [0, \infty)$
$$
w_n(t) \Df \Gamma\left(y + \int_0^{\cc}v_n(s) ds\right)(t).
$$From the Lipschitz property of $\Gamma$ (Condition
\ref{regular})
we have that
\[
\sup_{0 \le t < \infty}
|z_n(t) - w_n(t)| \le K|x-y|.
\]
Also clearly $w_n \in \acal(y)$.
Now let
$$
\tau'_n \Df \inf \{t \in (0, \infty): w_n(t) = 0\}.$$
Fix $n$ and suppose that $\tau_n \le \tau'_n$.
Then
\[
|w_n(\tau_n)| = |w_n(\tau_n) - z_n(\tau_n)|
 \le  K|x-y|.
\]
Hence from (\ref{bdd}), letting $C\doteq 4K^3\delta^{-1}$,
$$
\tau'_n \le \tau_n + C |x-y|.$$
Similarly it can be seen that if
$\tau'_n \le \tau_n$
then
$$
\tau_n \le \tau'_n + C |x-y|
$$
and thus
$$
|\tau_n - \tau'_n| \le C|x-y|.$$
We therefore have that
\[
\tau_n  \le  \tau'_n + C|x-y|
\le  T(y) + C|x-y|.
\]

Sending $n \to \infty$, it follows from (\ref{limit}) that
$$
T(x) \le T(y) + C|x-y|.$$
Since the role of $x$ and $y$ can be reversed, we have that
$$
|T(x) - T(y)| \le C|x-y|,
$$
and since $x$ and $y$ are arbitrary we have (i).

Next we show (ii).
Fix some $v \in \C(\delta)$,
and let $x \in G\backslash\{0\}$ be given.
With $\ii: [0, \infty) \to [0, \infty)$ denoting the
identity map,
clearly the trajectory $\{\Gamma(x + v\ii)(t)\}_{t \ge 0}$ belongs
to $\acal(x)$. Note that
\beq
\sup_{0 \le t \le M}
|\Gamma(x + v\ii)(t) - x|
&=& \sup_{0 \le t \le M}
|\Gamma(x + v\ii)(t) - \Gamma(x + 0\ii)(t)|\\
&\le& K M|v|.
\eeq
Therefore for any $M<|x|/K|v|$
\[
\inf_{0 \le t \le M}
|\Gamma(x + v\ii)(t)|
 \ge
|x| - KM|v|
> 0,
\]
which implies that $T(x) > M$.
Taking the supremum over $M<|x|/K|v|$ gives
\beq
T(x) & \ge & \frac{|x|}{K|v|}.
\eeq
This proves (ii) with $c=1/K|v|$.

Finally we prove (iii). Let $t>0$ be fixed.
 If $T(z(s)) = 0$ for some $s \in [0,t]$ then the result is
obviously true.  Now suppose that $T(z(s)) > 0$ for all
$s \in [0,t]$.
Let $\beta > 0$ be arbitrary and $u \in \acal(z(t))$ be such that
$\tau \Df \inf \{s \in [0, \infty) : u(s) = 0\}$ satisfies
$\tau > T(z(t)) - \beta$.
Define $\ti z: [0, \infty) \to \R^k$ by
\[
\ti z(s) =\left\{ \begin{array}{ll} z(s) &  s \le t\\
u(s-t) &  s > t.\end{array}
\right.
\]
Then $\ti z \in \acal(x)$ and
\beq
T(x) & \ge & \inf \{s \in [0, \infty): \ti z(s) = 0\}\\
& =& t + \tau \\
& \ge & T(z(t))+ t - \beta .
\eeq
Since $\beta > 0$ is arbitrary we have that
$T(z(t)) \le T(x) - t$.  This proves
the lemma.
\ink

\section{Stability of Constrained Diffusion Processes}
We begin with the following lemma.  For $x \in G$, let $\Om_0(x)$
be a $P-$null set such that for all $\om \not \in \Om_0(x)$
and $0 \le u < t < \infty$, $X^x(\cdot)=X^x(\cdot,\omega)$ satisfies
$$
X^x(t) = \Gamma \left(X^x(u) + \int_0^{\cdot}
  b(X^x(u+s)) ds
+ \int_0^{\cdot} \sigma (X^x(u+s)) dW_u(s) \right)(t-u),
$$
where $W_u(s) \Df W(s+u)$.
\begin{lemma}
\label{deltastep}
Let $T$ be the function defined in (\ref{hittime}).  Fix $x \in
G\setminus A$ and let
$\{X^x(t)\}_{ t \ge 0}$ be as in Theorem \ref{main}.
Let $\Delta > 0$ and $u > 0$ be arbitrary.
Fix $\om \not \in \Om_0(x)$.  Suppose that
$X^x(t, \om) \in G\backslash A$ for all $t \in (u, u + \D]$.
Then
\[
T(X^x(u + \D, \om)) \le (T(X^x(u, \om)) - \D)^+
+ KC \ov \nu(\om) ,
\]
where $C$ is as in Lemma \ref{hit} (i) and
\be
\label{nunov}
\ov \nu \Df \sup_{u \le s \le u +\D} \left |\int_{u}^s
\sigma(X^x(s)) dW(s)\right |.\ee
\end{lemma}
{\bf Proof:}
In the proof we will suppress $\om$ from the notation.
We begin by noticing that for $t \in [u , u + \D)$,
$X^x(t) = \ti X(t - u)$, where
\beq
\ti X(t) &\Df& \Gamma\left(X^x(u)
+ \int_0^{\cdot} b(X^x(s+ u))
ds \right . \\
& &\mbox{}+ \left . \int_0^{\cdot} \sigma(X^x(s+ u)) dW_u(s)
\right)(t); \;\; 0 \le t \le \Delta .
\eeq
Now define a sequence of $\R^k$ valued stochastic processes
$\{\tilde Y(t)\}_{0 \le t \le \D}$ as follows.
\be
\label{surro}
\tilde Y(t) \Df \Gamma \left( X^x(u) +
\int_0^{\cdot} b(X^x(s+u))
ds\right) (t).
\ee
Note that $\tilde Y(t)$ has absolutely continuous paths $P$-a.s.,
and that  $b(X^x(s+u))\in {\cal C}(\delta)$ for all $s \in [0,\Delta]$.
Also, note that by Condition \ref{regular} we have
\beqn
\sup_{0 \le t \le \D}
|\ti X(t) - \tilde Y(t)|
& \le &
K \sup_{0 \le t \le \D}
\left | \int_0^t \sigma(X^x(s+u)) dW_u(s)\right| \non \\
&=& K \ov \nu.
\label{lateruse}
\eeqn
Using the Lipschitz property of $T$ (Lemma \ref{hit} (i))
we have that
\beq
T(X^x(u + \D)) & = & T(\ti X(\Delta))\\
& \le & T(\tilde Y(\D)) + KC \ov \nu\\
& \le & (T(X^x(u)) - \D)^+ + KC \ov \nu,
\eeq
where the last inequality follows from Lemma \ref{hit}(iii).
\ink

\begin{lemma}
\label{doob}
Suppose that  $\{\alpha_i(t)\}$; $i = 1, 2, \cdots , l$ are
$\R^k$ valued $\sigma \{W(s): 0 \le s \le t
\}$-progressively measurable processes such that there exists
$\ov \alpha \in (0, \infty)$ for which
$$
|\alpha_i(t)| \le \ov \alpha ,$$
for all $t \in (0, \infty)$, $i \in \{1, \cdots , l\}$,
$P$-a.s.
Then for $\lambda \in (0, \infty)$
$$
\E\left( e^{\lambda  \sum_{i=1}^l
\left |\int_0^t \inn{\alpha_i(s)} {dW(s)} \right |}\right)
\le 2 e^{\frac{l^2\lambda^2 \ov \alpha^2 t}{2}},
$$
where $\inn{\cdot}{\cdot}$ is the usual inner product in $\R^k$.
\end{lemma}

\noindent
{\bf Proof:}
We first consider the case when $l =1$.
Observe that
$$E\left(\exp\left(\lambda \int_0^t \inn{\alpha_1(s)}{dW(s)} - \half \la^2 \int_0^t |\alpha_1(s)|^2 ds\right) \right) = 1$$
and
$$E\left(\exp\left(-\lambda \int_0^t \inn{\alpha_1(s)}{dW(s)} - \half \la^2 \int_0^t |\alpha_1(s)|^2 ds\right) \right) = 1.$$
Using the upper bound on $|\alpha_1(\cc)|$ we now have that
\beq
E\left(\exp\left(\lambda \left|\int_0^t \inn{\alpha_1(s)}{dW(s)}\right|\right) \right)
&\le & E\left(\exp\left(\lambda \int_0^t \inn{\alpha_1(s)}{dW(s)}\right) \right)\\
& +& E\left(\exp\left(-\lambda \int_0^t \inn{\alpha_1(s)}{dW(s)}\right) \right)\\
& \le & e^{\frac{\lambda^2 \ov \alpha^2 t}{2}} + e^{\frac{\lambda^2 \ov \alpha^2 t}{2}}\\
& = & 2 e^{\frac{\lambda^2 \ov \alpha^2 t}{2}} .
\eeq

This proves the lemma for the case $l=1$. Now we consider the case $l > 1$.
Note that
\beq
\E \left( e^{\lambda \sum_{i=1}^l |\int_0^t
    \inn{\alpha_i(s)}{dW(s)}|}\right)
& \le &
\E \left( \prod_{i=1}^l e^{\lambda |\int_0^t
    \inn{\alpha_i(s)}{dW(s)}|}\right)\\
& \le & \left( \prod_{i=1}^l \E\left(e^{l\lambda |\int_0^t
    \inn{\alpha_i(s)}{dW(s)}|}\right)
\right)^{\frac{1}{l}}\\
& \le & 2 \left( e^{\frac{l^3 \lambda^2 \ov{\alpha}^2
      t}{2}}\right)^{\frac{1}{l}}\\
&= & 2e^{\frac{l^2\lambda^2 \ov{\alpha}^2
      t}{2}}.
\eeq
\ink

In what follows, we will denote the set of positive integers by $\INT$.
\begin{lemma}
\label{doob-main}
Let $x \in G$ and $\D > 0$ be fixed.  For $n \in \INT$ let
$\nu_n$ be defined as follows:
\be
\label{nun}
\nu_n \Df \sup_{(n-1)\D \le s \le n \D} \left |\int_{(n-1)\D}^s
\sigma(X^x(s)) dW(s)\right |.\ee
Then for any
$\kappa \in (0, \infty)$ and $m, n \in \INT$; $m \le n$,
$$
\E\left( e^{\kappa \sum_{i=m}^n \nu_i}\right)
\le
\left [ 2\sqrt{2} e^{k^2\kappa^2 \gamma^2 \D} \right]^{(n-m+1)},
$$
where $\gamma$ is as in Condition \ref{growlip}.
\end{lemma}
{\bf Proof:}
For $t > 0$, let
$$
\G_t \Df \sigma \{ W(s): 0 \le s \le t \}.$$
Then
\be
\label{cond1}
\E e^{\kappa \sum_{i=m}^n \nu_i}
=
\E\left( e^{\kappa \sum_{i=m}^{n-1} \nu_i}
\left(
\E\left(e^{\kappa \nu_n} \mid \G_{(n-1)\D}\right)\right) \right).
\ee
Now
\beq
\E\left(e^{\kappa \nu_n} \mid \G_{(n-1)\D}\right)
& =&
\E \left(
\left. \sup_{(n-1)\D \le s \le n \D}
e^{\kappa \left|\int_{(n-1)\D}^s \sigma(X^x(u))dW(u)\right|}\right|
\G_{(n-1)\D} \right) \\
&=&
\E \left(
\left. \sup_{(n-1)\D \le s \le n \D}
e^{\kappa \left|\int_{(n-1)\D}^s \sigma(X^x(u))dW(u)\right|}\right|
X^x((n-1)\D) \right),
\eeq
where the last step follows from the Markov property of
$X^x$.

An application of Doob's maximal inequality for submartingales
yields that
the last expression is bounded above by
$$
2 \left( \E \left( \left. e^{2\kappa \left|\int_{(n-1)\D}^{n\D} \sigma(X(u))
      dW(u)\right |}
\right| X^x((n-1)\D) \right) \right)^{\half}.$$
By an application of Lemma \ref{doob}
and the observation that for positive real numbers $x_1, \cdots x_k$, $\sqrt{\sum_{i=1}^k x_i^2} \le
\sum_{i=1}^k x_i$ we have that the last expression
is bounded  above by
$$
2\sqrt{2}e^{k^2\kappa^2 \gamma^2 \D}.$$
Using this observation in (\ref{cond1}) we have the result by
iterating.
\ink

\skp
For $\D > 0$ let
\[
B^{\D} \Df \{y \in G: T(y) \le \D \}.
\]
Let $\{X^x(t)\}_{t \ge 0}$ be as in
Theorem \ref{uniqsoln}.  Given a compact set $B \subset G$, let
\be
\label{taub}
\tau_B(x) \Df \inf \{t: X^x(t) \in B \}.\ee
\begin{theorem}
\label{finitethm}
Let $\{X^x(t)\}_{t \ge 0}$ be as in
Theorem \ref{main}.
  Then there exists  $\Delta \in ( 0,\infty)$
such that for all $M \in (0, \infty)$
$$
\sup_{x: |x| \le M} \E(\tau_{B^{\D}}(x)) < \infty .$$
\end{theorem}

\noindent
{\bf Proof:}
Without loss of generality
we can assume that  $\Delta$ is chosen large enough so that
$B^{\Delta} \supset A$.

Let $\Delta \in (0, \infty)$ and let
$$A_n \Df \{\omega : \inf_{s \in [0, n\D]}T(X^x(s)) > \Delta \}. $$
Then
\beqn
P(A_n) & \le &
P\left( \Delta < T(X^x(n\D)) \le T(x) - n \D + CK\sum_{j=1}^n \nu_j
\right), \non \\
\label{an}
\eeqn
where $\{\nu_j\}_{j=1}^n$ are as in (\ref{nun}) and the
inequality follows from Lemma \ref{deltastep}.
Next observe that the probability on the right side of
(\ref{an})
is bounded above by:
\beq
P \left ( CK \sum_{i=1}^n \nu_i \ge (n+1)\D - T(x) \right)
& \le &
\frac{\E(e^{\alpha CK \sum_{i=1}^n \nu_i})}
{e^{\alpha ((n+1)\D - T(x))}} \\
& \le &
\frac{\left(2\sqrt{2}e^{k^2\alpha^2C^2K^2\gamma^2 \D}\right)^n}
{e^{\alpha ((n+1)\D - T(x))}}\\
& = &
\frac{e^{\alpha T(x)}}{e^{\alpha \D}}
e^{\left( k^2\alpha^2 C^2K^2 \gamma^2 - \alpha + \frac{\log
      8}{2\D}\right)
n \D},
\eeq
where $\alpha > 0$ is arbitrary
and the next to last inequality follows from Lemma
\ref{doob-main}.
Choose $\D > 0$ (sufficiently large) and $\alpha > 0$ (sufficiently
small)
 such that
$$
k^2 \alpha^2 C^2K^2 \gamma^2 - \alpha + \frac{ \log 8}{2\D} \Df -
\eta  < 0.$$
Then
\beq
P(X^x(s) \not \in B^{\D}; 0 \le s \le n\D )
& = & P(A_n)\\
& \le & \frac{e^{\alpha T(x)}}{e^{(\alpha - \eta)\D}}e^{-\eta (n+1)
  \D},
\eeq
for all $n \in \INT$.
Now let $t \in (0, \infty)$ be arbitrary and $n_0$ be such that
$t \in [n_0 \D , (n_0 + 1)\D]$.
Then
\beq
P(\tau_{B^{\D}}(x) > t)
& =&
P(X^x(s) \not \in B^{\D}; 0 \le s \le t) \\
& \le &
P(X^x(s) \not \in B^{\D}; 0 \le s \le n_0\D) \\
& \le &
\frac{e^{\alpha T(x)}}{e^{(\alpha - \eta)\D}}
e^{- \eta (n_0 + 1)\D} \\
& \le &
\frac{e^{\alpha T(x)}}{e^{(\alpha - \eta)\D}}
e^{- \eta t}.
\eeq
Hence
\beq
\E(\tau_{B^{\D}}(x)) &=&
\int_0^{\infty} P(\tau_{B^{\Delta}}(x) > t) dt \\
& \le & \frac{e^{\alpha T(x)}}{e^{(\alpha -
      \eta)\D}}\int_0^{\infty}e^{-\eta t} dt \\
& =& \frac{e^{\alpha T(x)}}{\eta e^{(\alpha - \eta) \D}}.
\eeq
Recalling that $T(\cc)$ is a continuous function we have from the
above inequality that for all $M \in (0, \infty)$
$$
\sup_{x: |x| \le M} \E(\tau_{B^{\D}}(x)) < \infty .
$$
\ink

\begin{lemma}
\label{tight}
For $x \in G$ let $\{X^x(t)\}_{t \ge 0}$ be as in
Theorem \ref{main}.  Then for all
$M \in (0, \infty)$ the family $\{X^x(t); t \ge 0, |x| \le M\}$
is tight.
\end{lemma}
{\bf Proof:}
Let $\Delta > 0$ be large enough so that $B^{\Delta} \supset A$.
Fix $\om \in \Om_0(x)$, where $\Om_0(x)$ is as defined at the
beginning of this section.
In the rest of the proof we will suppress the dependence of all
random variables on $\om$ in the notation.
Let
$$S(\Delta) \Df \{j \in \{1, 2, \cdots , n-1\}: T(X^x(t))\le
\D\;\;
\mbox{for some} \; t \in [(j-1)\Delta , j \Delta)\}.$$
Define
\[
m = \left\{ \begin{array}{ll}\max \{j: j \in S(\Delta)\} &  \mbox{if $S(\D)$ is
  non empty}\\
 0 & \mbox{otherwise}.
\end{array}
\right.
 \]
 From Lemma \ref{deltastep} we
have that
\be T(X^x(n\D)) \le  T(X^x(m\Delta)) +
\sum_{j=m+1}^n(KC\nu_j - \D).
\label{z1}
\ee
Let
$$t \Df \sup \{s \in [(m-1)\D , m \D): T(X^x(s)) \le \D \}.$$
If $m > 0$, we have from  Lemma \ref{deltastep} that
\beq
T(X^x(m\Delta))
& \le &
(T(X^x(t)) - (m\D - t))^+ +
KC \sup_{t < s < m\Delta}\left| \int_t^s \sigma(X^x(u)) dW(u)
\right|\\
& \le &
 \D + KC \sup_{t < s < m\Delta}\left| \int_t^s \sigma(X^x(u)) dW(u)
\right|\\
& \le & \D + 2 KC \sup_{(m-1)\Delta < s < m\Delta}\left|
\int_{(m-1)\D}^s \sigma(X^x(u)) dW(u)
\right|\\
&=& \D + 2KC\nu_m,
\eeq
where in obtaining the second inequality we have
used the fact that for $(m-1)\D \le t \le s \le m \D$
$$\left|\int_t^{s}\sigma(X^x(u)) dW(u)\right|
    \le \left |\int_{(m-1)\D}^{s}\sigma(X^x(u)) dW(u)\right|
+ \left |\int_{(m-1)\D}^{t} \sigma(X^x(u)) dW(u)\right|.$$

Using this observation in (\ref{z1}) we have that
\beq
T(X^x(n\D)) & \le & T(x) + 2\D + \sum_{j=m}^n (2KC\nu_j(x) - \D)\\
& \le & T(x) + 2\D + \max_{1 \le
l \le n} \sum_{j=l}^n (2KC\nu_j(x) - \D), \eeq
where we have
written $\nu_j \equiv \nu_j(x)$ in order to explicitly bring out
its dependence on $x$. Hence for $M_0 \in (0, \infty)$ \beq
P(T(X^x(n\D)) \ge M_0) &\le& P\left( \max_{1 \le l \le n}
\sum_{j=l}^n(2KC\nu_j(x) - \D) \ge M_0 - T(x) - 2\D\right)\\ &\le &
\sum_{l=1}^n P\left(2KC\sum_{j=l}^n \nu_j(x) \ge M_0 + (n-l-1)\D -T(x)\right)\\
& \le & \frac{e^{\alpha (T(x) + \D)}}{e^{\alpha M_0}}\sum_{l=1}^n
\frac{\E(e^{\alpha 2KC \sum_{j=l}^n \nu_j(x)})}{e^{\alpha (n-l)
    \D}},
\eeq
where $\alpha > 0$ is arbitrary.
 From Lemma \ref{doob-main} we now have that
\beq
P(T(X^x(n\D)) \ge M_0)
& \le & \frac{e^{\alpha (T(x) + \D)}}{e^{\alpha M_0}}\sum_{l=1}^n
\frac{\left(2\sqrt{2} e^{8k\alpha^2C^2K^2 \gamma^2
      \D}\right)^{n-l+1}}
{e^{\alpha (n-l)
    \D}},\\
& \le & \frac{e^{\alpha (2\D + T(x))}}{e^{\alpha M_0}}\sum_{l=1}^n
 e^{(8k\alpha^2C^2K^2 \gamma^2 \D - \alpha \D + \frac{\log
     8}{2})(n-l +1)}.
\eeq
Now choose $\alpha$ and $\D$ so that
$$
8k\alpha^2C^2K^2 \gamma^2 \D - \alpha \D + \frac{\log
     8}{2}  = - \theta < 0.$$
Then
\beq
P(T(X^x(n\D)) \ge M_0) &\le & \frac{e^{\alpha (2\D + T(x))}}{e^{\alpha M_0}}
\sum_{l=1}^n
e^{-\theta(n-l)}  \\
& \le & \frac{e^{\alpha (2\D + T(x))}}{e^{\alpha M_0}(1- e^{-\theta})}.
\eeq
Hence for all $M, M_0 \in (0, \infty)$
$$
\sup_{n \in \INT , |x| \le M}
 P(T(X^x(n\D)) \ge M_0) \le \frac{e^{\alpha (2\D + \frac{4K^2}{\delta }M)}}
 {e^{\alpha
    M_0}(1-e^{-\theta})}.
$$
 From Lemma \ref{hit} (ii) it now follows that
\be
\label{onebd}
\sup_{n \in \INT, |x| \le M} P(X^x(n\D) \ge M_0) \le
\frac{e^{\alpha (2\D + \frac{4K^2}{\delta }M)} }{e^{\alpha
    c M_0}(1-e^{-\theta})}.
\ee
Now let $t \in [n \D, (n+1)\D]$ and consider the process
$\{\tilde Y(t)\}_{0 \le t \le \D}$
defined in
(\ref{surro}) with $u$ there replaced by $n\D$.
For each $n$ it follows from (\ref{lateruse}) that
$$
|X^x(t) - \tilde Y(t - n \D)| \le K \nu_n(x).$$
Define a function $\tilde b: G \to \R^k$ which agrees
with $b$ off $A$ and satisfies Condition \ref{permit} with $A=\emptyset$
as well as equations (\ref{lipbsig}).
Also define
\[
 Y^*(t) \Df \Gamma \left( X^x(n\D) +
\int_0^{\cdot} \tilde b(X^x(s+ n\D))
ds\right) (t).
\]
Clearly
$$ L \Df \sup_{x \in G}|b(x) - \tilde b(x)| < \infty .$$
Furthermore
\beq
 |\tilde Y(t - n\D)| & \le &
|\tilde Y(t - n\D) - Y^*(t - n\D)| + |Y^*(t - n\D)|\\
& \le &
KL\D + K |X^x(n\D)|,
\eeq
where in obtaining the last inequality we have used
the Lipschitz property of $\Gamma$ and Theorem \ref{staba}.
Combining the above observations we have that
$$
|X^x(t)| \le K( \nu_n(x) + |X^{x}(n\D)|) + KL\D.$$
Therefore for $M_0 \in (0, \infty)$, any $n$, and $t \in [n\D, n\D+\D]$
\beqn
P(|X^x(t)| \ge M_0) \le P\left(\nu_n(x) \ge \frac{M_0 - KL\D}{2K}\right)
+ P\left(|X^x(n\D)| \ge \frac{M_0- KL\D}{2K}\right).\non \\
\label{secbd}
\eeqn
Clearly, the family $\{|\nu_n(x)|; n \ge 1, |x| \le M \}$ is tight.
Now let $\eta > 0$ be arbitrary.  Choose $M_0 \in (0, \infty)$
such that
$$\sup_{n \in \INT , |x| \le M} P\left(|\nu_n(x)| \ge \frac{M_0 -KL\Delta }{2K}\right)
 \le \frac{\eta}{2}$$
and
$$\frac{e^{\alpha(2\D + \frac{4K^2}{\delta }M)}}
{e^{\frac{\alpha c(M_0 -KL\D) }{2K}}(1-e^{-\theta})} \le
\frac{\eta}{2}.$$
Then from (\ref{onebd}) and (\ref{secbd}) we have that
$$
\sup_{t \ge 0, |x| \le M} P(|X^x(t)| \ge M_0) \le \eta .$$
Since $\eta > 0$ is arbitrary, we have the result.
\ink

\skp
 From Theorem \ref{finitethm}, Lemma \ref{tight} and Condition \ref{irred}
the proof of positive recurrence and the existence and uniqueness
of an invariant measure for $\{X^x(t)\}_{t \ge 0}$ is standard
(cf. \cite{dupwil}).  However for the sake of completeness we
present the proof below.

\skp \noi
{\bf Proof of Theorem \ref{main}:}
 Denote the measure
induced by $\{X^x(\cc)\}$ on $C([0,\infty): G)$ by $P_x$, where $C([0,\infty): G)$ is
the space of $G$ valued continuous functions defined on the nonnegative real line.
In arguments that are presented below,
it will be convenient to let initial conditions be defined through conditioning,
rather than though the superscript as in $X^x$.
As a consequence,
instead of $X^x$ we will work with the canonical process $\xi(\cc)$ on $C([0,\infty): G)$, and the canonical
filtration which we denote by $\{\F_t\}$.
Finally, the expectation operator corresponding to the probability measure
$P_x$ will be denoted by $\E_x$.
Given a compact set $B \subset G$, let
\[
\ti \tau_B \Df \inf \{t: \xi(t) \in B \}.\]
In order to show positive recurrence, we need
to
show
 that if $S$ is an arbitrary compact set in $G$ with
positive Lebesgue measure then for all $x \in G$,
$\E_x \ti \tau_S < \infty$.
Let $B^{\D}$ be as in Theorem \ref{finitethm},
and
let $r \in (0, \infty)$ be such that $B^{\D} \subset \{x: |x| \le
r\}$.
Then from Theorem \ref{finitethm} we have that for all
$C \in (0, \infty)$
\be \label{ball}
\sup_{x: |x| \le C}\E_x(\ti \tau_{B_r}) < \infty ,
\ee
where $B_r \Df \{x \in G: |x| \le r\}$.
 From the uniform non degeneracy assumption (Condition
\ref{irred}),
we have (cf. \cite{harwil})
$$  p(S) \doteq \inf_{x \in B_r} P_x(\xi(1) \in S) > 0.$$
Furthermore, Feller property of
$\{X^x(\cdot)\}$ implies that the family
$\{X^x(t): x \in B_r, 0 \leq t \leq 1\}$
is tight, and
 so there
exists $M \in (0, \infty)$ such that
\be
\label{apfin1}
\inf_{x \in B_r}
P_x(\xi(1) \in S \;\; \mbox{and}\;\;
|\xi(t)| \le M \;\; \mbox{for all}\;\; t \in [0,1])
\ge \frac{p(S)}{2}.
\ee
Let $C \in (M, \infty)$ be fixed, and define
\[ \hat \tau \Df \inf\{t: |\xi(t)| \ge C\} \]
and $\ti \tau \Df \min \{1, \hat \tau, \ti \tau_S \}$.
If $y \in B_r$, then by the strong Markov property
\beqn
\E_y(\ti \tau_S) & = &
\E_y \left ( \E_y \left ( \ti \tau_S \mid \F_{\ti \tau}\right) \right)\non \\
&\le& \E_y \left( \ti \tau + \E_{\xi(\ti \tau)}(\ti \tau_S) \right)\non \\
& \le & 1 + \E_y \left(\E_{\xi(\ti \tau)}(\ti \tau_S) \right). \label{apfin2}
\eeqn
Now define
\[ \Lambda \Df \{ \xi(\cc) \in C([0, \infty): G): \sup_{0 \le t \le 1}|\xi(t)| \le M \;\mbox{and}\; \xi(1) \in S \}.\]
Since $E_{\xi(\tilde \tau)}(\tilde \tau_S)=0$ w.p.1 if $\xi \in \Lambda$,
for $y \in B_r$
\beqn
\E_y\left(\E_{\xi(\ti \tau)}(\ti \tau_S) \right) &=&
\E_y\left(\I_{\Lambda}(\xi)\E_{\xi(\ti \tau)}(\ti \tau_S) \right) +
\E_y\left(\I_{\Lambda^c}(\xi)\E_{\xi(\ti \tau)}(\ti \tau_S) \right) \non \\
&=& \E_y\left(\I_{\Lambda^c}(\xi)\E_{\xi(\ti \tau)}(\ti \tau_S) \right). \label{apfin3}
\eeqn
Next, fix $z \in G$ such that $|z| \le C$.  Then
\[ \E_z(\ti \tau_S) = \E_z(\ti \tau_{B_r} + (\tau_S - \ti \tau_{B_r})). \]
Once more using the strong Markov property, we have
\beqn
\E_z(\ti \tau_S) & \le & \E_z\left(\ti \tau_{B_r} + \sup_{x \in B_r} \E_x(\ti \tau_S)\right)\non \\
& \le & \sup_{z: |z| \le C} \E_z(\ti \tau_{B_r}) + \sup_{x \in  B_r} \E_x(\ti \tau_S).
\label{apfin4}
\eeqn
Observing that $|\xi(\ti \tau)| \le C$ and combining (\ref{apfin3}) and (\ref{apfin4}) we have that
for $y \in B_r$
\be
\E_y\left(\E_{\xi(\ti \tau)}(\ti \tau_S) \right)
\le \left( \sup_{z: |z| \le C} \E_z(\ti \tau_{B_r}) + \sup_{x \in B_r} \E_x(\ti \tau_S)\right) P_y(\Lambda^c).
\label{apfin5}
\ee
From (\ref{apfin1}), (\ref{apfin2}) and (\ref{apfin5}) it now follows that
\[
\sup_{y \in B_r} \E_y(\ti \tau_S)
\le  1 + \sup_{z: |z| \le C} \E_z(\ti \tau_{B_r}) + \left(1 - \frac{p(S)}{2}\right)\sup_{x \in B_r}\E_x(\ti \tau_S).
\]
Thus
\[
\sup_{y \in B_r} \E_y(\ti \tau_S) \le \frac{2}{p(S)}\left\{ 1 + \sup_{z: |z| \le C} \E_z(\ti \tau_{B_r}) \right \}.
\]
A final application of the strong Markov property now yields that for $x \in G$
\beq
\E_x(\ti \tau_S) & \le &
\E_x(\ti \tau_{B_r}) +  \sup_{y \in B_r} \E_y(\ti \tau_S) \\
& \le &
\E_x(\ti \tau_{B_r}) + \frac{2}{p(S)}\left\{ 1 + \sup_{z: |z| \le C} \E_z(\ti \tau_{B_r}) \right \}\\
&< & \infty ,
\eeq
where the last inequality follows from Theorem \ref{finitethm} and (\ref{ball}).
This completes the proof of positive recurrence.

Finally we consider the existence and uniqueness of invariant
measures.
 From Lemma \ref{tight} we have that the family of measures
$\{\mu_t ; t \ge 1\}$ defined by
$$\mu_t(B) \Df \frac{1}{t}\int_0^t P(X^x(s) \in B) ds$$
is tight.  Since the Markov process $\{X^x(t)\}$ is Feller we
have that any weak limit of $\{\mu_t\}$ is an invariant measure.
(See, for example,  the proof of Theorem 4.1.21, Chapter I, \cite{skor}).
Finally uniqueness follows as in \cite{harwil, dupwil} in view of
Condition \ref{irred}.
\ink

\begin{remark}
\label{finrem}
The Lipschitz and growth condition (Condition \ref{growlip}) on
$b$
and $\sigma$ are essentially assumed to guarantee a unique
solution
to the constrained diffusion process (\ref{maineq}) which is
Feller Markov.
The conclusion of
Theorem
\ref{main} continues to hold with  the same proof if
Condition \ref{growlip} is replaced by the assumption that
(\ref{growsig})
holds for some $\gamma \in (0, \infty)$, $b$ is locally bounded and
(\ref{maineq}) has a unique weak solution with continuous paths
for every $x \in G$
and the solution is Feller-Markov.
\end{remark}
\section*{Acknowledgments.} We will like to thank the referees for a careful review of the paper.


\Line{\AOPaddress{Department of\\ Electrical Engineering,\\
Technion,\\Israel Institute\\ of Technology,\\
Haifa 32000, Israel}\hfill
\AOPaddress{Department of\\ Statistics,\\
University of\\ North Carolina\\ at Chapel Hill,\\
NC 27599-3260}\hfill
\AOPaddress{Lefschetz Center\\ for Dynamical Systems,\\
Division of\\ Applied Mathematics,\\
Brown University,\\
Providence, RI 02912}
}

\end{document}